\documentclass[10pt]{article}
\usepackage[ansinew]{inputenc}
\usepackage[francais]{babel}
\usepackage{amssymb,hyperref}
\usepackage{graphicx}
\hypersetup{pdfpagemode=FullScreen,colorlinks=true}

\title{Plongements quasiisométriques du groupe de Heisenberg dans $L^p$, d'après Cheeger, Kleiner, Lee, Naor}
\author{Pierre Pansu$^{1,2}$\footnote{$^{1}$ Univ Paris-Sud, Laboratoire de Mathématiques d'Orsay, Orsay, F-91405 ;\hfill\eject\indent\hskip7pt $^{2}$ CNRS, Orsay, F-91405.}}

\date{\today}

\newtheorem{thm}{Théorème}
\newtheorem{lem}{Lemme}
\newtheorem{prop}[lem]{Proposition}
\newtheorem{cor}[lem]{Corollaire}
\newtheorem{defi}[lem]{Définition}
\newtheorem{rem}[lem]{Remarque}
\newtheorem{nota}[lem]{Notation}
\newtheorem{ques}[lem]{Question}
\newtheorem{prob}[lem]{Problème}

\newtheorem{conc}[lem]{Conclusion}
\newtheorem{ex}[lem]{Exemple}
\newtheorem{exs}[lem]{Exemples}

\def\proof{\par\medskip\noindent {\bf Preuve}}

\def\R{\mathbb{R}}

\def\Z{\mathbb{Z}}
\def\H{\mathbb{H}}

\def\n#1{{|\hskip-1pt | #1 |\hskip-1pt |}}
\def\sp{{\rm sp}\,}

\def\texorpdfstring#1#2{#1}

\begin{document}

\maketitle

\begin{quote}
{\small
ABSTRACT. This is a short survey of Cheeger and Kleiner's nonembeddability theorem for Heisenberg group into $L^1$.
}
\end{quote}

\section{Introduction}

Cet exposé porte sur une question de géométrie métrique. Etant donné un espace métrique $(X,d)$, à quelle condition peut on plonger de façon bilipschitzienne l'espace $(X,d^{\alpha})$, $0<\alpha<1$, dans un espace modèle, par exemple, un espace de Banach d'usage courant ? 

\subsection{Résultats}

D'une certaine façon, la question a été résolue par P. Assouad.

\begin{defi}
Un espace métrique est {\em doublant} s'il existe une mesure $\mu$ et une constante $C$ telles que pour toute boule $B$, $\mu(2B)\leq C\,\mu(B)$.
\end{defi}

\begin{ex}
Le groupe d'Heisenberg $(\H,d)$ muni de sa métrique de Carnot est doublant.
\end{ex}

\begin{thm}
{\em (P. Assouad, 1982, \cite{Assouad-doublant})}. Soit $(X,d)$ un espace métrique doublant. Alors pour tout $\alpha<1$ et tout $p\geq 1$, $(X,d^{\alpha})$ admet un plongement bilipschitzien dans $L^p$.
\end{thm}
En fait, Assouad fournit un plongement en dimension finie. La constante de Lipschitz dépend de $C$, de $\alpha$ et de $p$. On va voir que cette dépendance est intéressante. On ne peut pas impunément faire tendre $\alpha$ vers 1. 

Elever une distance à une puissance $<1$ fait disparaître toutes les courbes de longueur finie. On peut penser que celles-ci constituent un obstacle au plongement bilipschitzien. Après l'espace euclidien (qui se plonge évidemment), le premier exemple d'espace doublant possédant des courbes de longueur finie auquel on pense est le groupe d'Heisenberg.

\begin{thm}
\label{CK}
{\em (J. Cheeger, B. Kleiner, 2006, \cite{Cheeger-Kleiner-L1})}. Le groupe d'Heisenberg $(\H,d)$ n'admet pas de plongement bilipschitzien dans $L^1$.
\end{thm}

Ce résultat répond à une question posée par A. Naor, elle-même motivée par une conjecture venant de l'informatique théorique. Comme on le verra plus loin, pour l'informatique, le théorème de non plongeabilité prend tout son sens quand on le confronte au résultat suivant, qui précise, dans le cas du groupe d'Heisenberg, le théorème d'Assouad.

\begin{thm}
{\em (J. Lee, A. Naor, 2006, \cite{Lee-Naor})}. Il existe une métrique invariante $d'$ sur le groupe d'Heisenberg, équivalente à $d$, telle que $(\H,d'^{1/2})$ se plonge isométrique\-ment dans $L^2$.
\end{thm}

\subsection{Organisation de l'exposé}

L'objet principal de ce laïus est de raconter la preuve du théorème \ref{CK}. Si le fait que le groupe d'Heisenberg ne se plonge pas dans $L^p$, $p>1$, résulte de mécanismes classiques, il me semble que le cas $p=1$ introduit des idées nouvelles. Dans la foulée, je présente des problématiques actuelles en informatique (je remercie A. Naor de me l'avoir patiemment expliqué) et en théorie des groupes.

\medskip
\textbf{Plan}.

\begin{enumerate}
  \item Introduction
  \item Motivation provenant de la théorie géométrique des groupes
  \item Motivation provenant de l'informatique théorique
  \item Résultats antérieurs sur le problème de plongement
  \item Preuve du théorème de non plongement dans $L^1$
\end{enumerate}

\section{Motivation provenant de la théorie géo\-mé\-tri\-que des groupes}

La théorie des groupes apporte des exemples (les groupes) d'espaces métriques qui ne sont pas doublants, pour la plupart. Néanmoins, on va voir que les résultats ci-dessus apportent un éclairage utile.

\subsection{Compression}

\begin{defi}
Une application $f:X\to Y$ entre espaces métriques est un {\em plongement uniforme} s'il existe une constante $C$ et une fonction $\phi$ tendant vers $+\infty$ telles que pour tous $x$, $x'\in X$,
\begin{eqnarray*}
\phi(d(x,x'))\leq d(f(x),f(x'))\leq C\,d(x,x').
\end{eqnarray*}
\end{defi}

\begin{thm}
{\em (G. Yu, 2000, \cite{Yu})}. La conjecture de Baum-Connes grossière est vraie pour les groupes qui possèdent un plongement uniforme dans un espace de Hilbert.
\end{thm}

\begin{thm}
{\em (M. Gromov, 2003, \cite{Gromov-RWRG})}. Il existe des groupes de présentation finie qui ne possèdent aucun plongement uniforme dans un espace de Hilbert.
\end{thm}

Ceci motive l'étude détaillée des plongements uniformes et la terminologie suivante

\begin{defi}
({\em M. Gromov, \cite{Gromov-AIIG})}. Le plus grand $\phi$ qui convient dans la définition d'un plongement uniforme $f$ est appelé la {\em compression} de $f$. 
\end{defi}

\begin{ques}
Etant donnés un espace métrique $X$ et une classe d'espaces métri\-ques $\mathcal{Y}$, quelle est la plus grande compression possible pour des plongements de $X$ dans des éléments de la classe $\mathcal{Y}$ ?
\end{ques}

\subsection{Exemples de plongements dans des espaces \texorpdfstring{$\ell^p$}{}}

\begin{thm}
{\em (J. Bourgain, 1986, \cite{Bourgain-86}, R. Tessera, 2006, \cite{Tessera})}. Supposons qu'il existe un plongement d'un arbre régulier dans $\ell^p$ de compression $\phi$. Alors 
$$
\displaystyle \int_{1}^{+\infty}(\frac{\phi(t)}{t})^{\max\{2,p\}}\frac{dt}{t}<+\infty.
$$

Réciproquement, toute fonction croissante satisfaisant 
$$
\displaystyle \int_{1}^{+\infty}(\frac{\phi(t)}{t})^{p}\frac{dt}{t}<+\infty
$$ 
est majorée par la compression d'un plongement dans $\ell^p$ d'un espace métrique quelconque tiré de la liste suivante : arbres réguliers, réseaux uniformes des groupes de Lie connexes, groupes hyperboliques, produits en couronne $F\wr\Z$ où $F$ est un groupe fini.
\end{thm}

\begin{ex}
{\em (Folklore)}. Les arbres se plongent dans $\ell^p$ avec compression $\phi(t)\geq t^{1/p}$.
\end{ex}
En effet, fixons une origine $o\in T^0$ et plongeons l'arbre $T$ dans $\ell^1 (T^0)$ en envoyant un sommet $x$ sur la fonction caractéristique de la géodésique de $o$ à $x$. Ensuite, on plonge $\ell^1 (T^0)$ dans $\ell^p (T^0)$ de la façon évidente. On obtient une compression $t^{1/p}$.

\begin{ques}
Est-ce que tous les plongements intéressants d'espaces métriques dans $\ell^p$ sont obtenus de cette manière, i.e. en passant par $\ell^1$ ?
\end{ques}

Pendant un certain temps, on aurait pu penser que oui. Par exemple, la {\em Propriété A} de Yu, qui donne une condition suffisante pour qu'un espace métrique se plonge uniformément dans $L^{2}$, fonctionne de cette façon. Les résultats combinés de Cheeger, Kleiner, Lee et Naor montrent que ce n'est pas le cas. On peut bien se plonger dans $L^2$ (après remplacement de $d$ par $\sqrt{d}$) sans se plonger dans $L^1$.

\section{Motivation provenant de l'informatique théo\-ri\-que}

Un plongement d'un espace métrique fini dans l'espace euclidien ou/et dans $\ell^1$ intervient souvent dans les algorithmes de base de l'informatique théorique. On va le voir pour le problème \texttt{Sparsest Cut}.

\subsection{La conjecture de Goemans et Linial}

Le meilleur algorithme connu pour une résolution approchée de \texttt{Sparsest Cut}, \texttt{SDP} (je vais l'expliquer), donne, dans le pire des cas, pour un graphe à $n$ sommets, une réponse qui diffère de l'optimum d'un facteur au plus égal à la constante $L_n$ définie comme suit.

\begin{nota}
On considère tous les espaces métriques à $n$ points $(X,d)$ tels que $(X,d^{1/2})$ se plonge isométriquement dans l'espace euclidien. Soit $L_n$ le plus petit $L$ tel que tous ces espaces admettent un plongement $L$-lipschitzien et qui augmente les distances dans $\ell^1$.
\end{nota}

M.X. Goemans (en 1997, \cite{Goemans}) et N. Linial (en 2002, \cite{Linial}) ont demandé si $L_n$ est borné indépendamment de $n$.

\medskip

En 2005, S. Khot et N. Vishnoi ont donné un contre-exemple, \cite{Khot-Vishnoi}. 

\medskip

Les boules du graphe de Cayley du groupe de Heisenberg discret fournissent d'autres contre-exemples, plus naturels et ayant des propriétés supplémentaires.

\subsection{\texorpdfstring{\texttt{Sparsest Cut}}{}}

\begin{prob}
\texttt{Sparsest Cut} consiste à calculer la constante de Cheeger d'un graphe fini.
\end{prob}
Une {\em coupure} dans un graphe pondéré $G$ est une partition des sommets en $G^0 =S\cup \bar{S}$.
\begin{eqnarray*}
\Phi(S)=\frac{\#\partial S}{\# S\,\#\bar{S}}.
\end{eqnarray*}
La {\em constante de Cheeger} de $G$ est $\displaystyle \Phi^* (G)=\min_{\emptyset\subsetneq S\subsetneq G^0}\Phi(S)$.

\medskip

Le calcul exact de $\Phi^*$ est NP-difficile. Mais une coupure (presque) optimale est fréquemment utilisée dans des algorithmes.

\medskip

\textbf{Arithmétisation}
\begin{eqnarray*}
\Phi(S)=\frac{\sum_{\mathrm{aretes}\,uv}m(uv)|1_S (u)-1_S (v)|}{\sum_{u}\sum_{v}|1_S (u)-1_S (v)|}.
\end{eqnarray*}
Soit $d(u,v)=|1_S (u)-1_S (v)|$. C'est une semi-distance sur $G^0$, induite par une application vers l'espace métrique à 2 points $\{0,1\}$. L'enveloppe convexe de ces semi-distances est exactement l'ensemble des semi-distances plongeables dans $\ell^1$. 
Par conséquent,
\begin{eqnarray*}
\Phi^* &=&\min_{d\,\mathrm{plongeable\, dans}\,\ell^1}\frac{\sum_{\mathrm{aretes}\,uv}m(uv)d(u,v)}{\sum_{u}\sum_{v}d(u,v)}\\
&=&\min\{\sum_{\mathrm{aretes}\,uv}m(uv)d(u,v)\,|\,d\,\mathrm{plongeable\, dans}\,\ell^{1},\,\sum_{u}\sum_{v}d(u,v)=1\}.
\end{eqnarray*}

\subsection{Approche de \texttt{Sparsest Cut} par la programmation liné\-aire}

Malheureusement, le problème de décider si un espace métrique fini est plongeable ou non dans $\ell^1$ est NP-complet.

\textbf{Relaxation}. Oublions la condition de plongeabilité dans $\ell^1$. Cela ramène à un problème de programmation linéaire, noté \texttt{LP}, pour lequel il existe des algorithmes polynômiaux. Soit $\Phi^{LP}$ le minimum de ce problème.

\begin{thm}
{\em (J. Bourgain, 1985, \cite{Bourgain-85})}. Tout espace métrique à $n$ points se plonge dans $L^2$ (et donc dans $L^1$) avec distorsion au plus $O(\log(n))$. C'est optimal.
\end{thm}

\begin{cor}
$$\Phi^{LP}\leq\Phi^{*}\leq C\,\log(n)\Phi^{LP},$$
ce qui montre que \texttt{LP} fournit une approximation de $\Phi^*$ à un facteur multiplicatif $\log(n)$ près.
\end{cor}
\proof. La métrique $d'$ plongeable dans $\ell^1$ qui est $O(\log(n))$-proche de la solution $d$ de \texttt{LP} satisfait 
$$\Phi^{*}\leq\Phi(d')\leq C\,\log(n)\Phi^{*} = C\,\log(n)\Phi^{LP}.$$

\textbf{Procédure d'arrondi}. {\em (N. Linial, E. London, Y. Rabinovich, 1995, \cite{Linial-London-Rabinovich})}. Le plongement de Bourgain est calculable en temps polynômial, on peut en tirer une coupure $S$ qui réalise approximativement $\Phi^{LP}$.


\subsection{Approche de \texorpdfstring{\texttt{Sparsest Cut}}{} via la programmation semi-définie}

\textbf{Arithmétisation}.
\begin{eqnarray*}
\Phi(S)=\frac{\sum_{\mathrm{aretes}\,uv}m(uv)|1_S (u)-1_S (v)|^2}{\sum_{u}\sum_{v}|1_S (u)-1_S (v)|^2}.
\end{eqnarray*}

\textbf{Relaxation}. On autorise des fonctions $x:G^0 \to L^2$ et non seulement à valeurs dans $\{0,1\}$, en gardant la contrainte
$$\forall u,\,v,\,w\in G^0 ,\quad |x(u)-x(v)|^2 \leq |x(u)-x(w)|^2 + |x(w)-x(v)|^2 ,$$
satisfaite par les fonctions caractéristiques. Cela ramène à un problème de programmation semi-définie, noté \texttt{SDP}, pour lequel il existe des algorithmes polynômiaux. Soit $\Phi^{SDP}=\Phi(x)$ son minimum.

Soit $d(u,v)=|x(u)-x(v)|^2$. C'est une semi-distance sur $G^0$, et $d^{1/2}$ est induite par un plongement dans l'espace euclidien. Par conséquent,
\begin{eqnarray*}
\Phi^{SDP}=\min\{\sum_{\mathrm{aretes}\,uv}m(uv)d(u,v)&|&d\,\mathrm{distance},\,\sqrt{d}\,\mathrm{plongeable\, dans}\,L^{2},\\
&&\sum_{u}\sum_{v}d(u,v)=1\}.
\end{eqnarray*}
Clairement,
\begin{eqnarray*}
\Phi^{SDP}\leq\Phi^{*}\leq L_n \Phi^{SDP},
\end{eqnarray*}
ce qui montre que $\Phi^*$ est calculable en temps polynômial à un facteur multiplicatif $L_n$ près.

\subsection{Estimation de \texorpdfstring{$L_n$}{}}

\begin{thm}
{\em (S. Arora, J. Lee, A. Naor, 2005, \cite{Arora-Lee-Naor})}. Soit $(X,d)$ un espace métrique à $n$ points. On suppose que $(X,d^{1/2})$ se plonge isométriquement dans $L^2$. Alors $(X,d)$ se plonge aussi dans $L^2$ avec distorsion $O(\sqrt{\log(n)}\log(\log(n)))$.
\end{thm}

\begin{rem}
C'est presque optimal, puisque l'ensemble des sommets du $n$-cube $\ell^1$ ne se plonge pas dans $L^2$ avec distorsion $<\sqrt{n}$ {\em (Enflo, 1969, \cite{Enflo})}.
\end{rem}

\begin{cor}
$L_n = O(\sqrt{\log(n)}\log(\log(n)))$.
\end{cor}
En effet, $L^2$ se plonge isométriquement dans $L^1$.

\begin{rem}
La non plongeabilité du groupe de Heisenberg entraîne une borne inférieure sur $L_n$. Cheeger, Kleiner and Naor annoncent qu'elle peut être rendue effective. Conjecturalement, $L_n =\Omega(\log(\log(n))^{\delta})$ pour  un $\delta>0$.
\end{rem}

\begin{conc}
L'approche SDP donne actuellement la meilleure solution con\-nue du problème \texttt{Sparsest Cut} général. Dans le cas particulier où les poids sont tous égaux, {\em S. Arora, E. Hazan, S. Kale, 2004, \cite{Arora-Hazan-Kale}} donnent un algorithme polynômial différent qui calcule $\Phi^*$ à un facteur $O(\sqrt{\log(n)})$ près.
\end{conc}

\section{Résultats antérieurs sur la plongeabilité quasiisométrique}
\subsection{Plongements dans \texorpdfstring{$\ell^p$}{} et \texorpdfstring{$L^p$}{}, \texorpdfstring{$p>1$}{}}

\begin{thm}
{\em (Semmes, 1996, \cite{Semmes})}. $\H$ ne se plonge pas quasiisométriquement dans un espace de Banach de dimension finie.

{\em (Pauls, 2001, \cite{Pauls})}. $\H$ ne se plonge pas quasiisométriquement dans un espace de Hilbert, ni, plus généralement, dans un espace $CAT(0)$.

{\em (Cheeger-Kleiner, 2006, \cite{Cheeger-Kleiner-L1})}. $\H$ ne se plonge pas quasiisométriquement dans un espace de Banach qui possède la propriété de Radon-Nikodym.
\end{thm}

\begin{defi}
{\em (Heinonen-Koskela, 1996, \cite{Heinonen-Koskela})}. On dit qu'un espace métrique est {\em PI} s'il est doublant et s'il satisfait une inégalité de Poincaré $(1,p)$,
\begin{eqnarray*}
\oint_B |f-\oint_{B}f|\leq\mathrm{const.~diameter}(B)\left(\oint_{2B}|\nabla f|^{p}\right)^{1/p},
\end{eqnarray*}
pour tous les sur-gradients $|\nabla f|$ de $f$.
\end{defi}

\subsection{Propriété GFDA}

\begin{defi}
Un espace de Banach $V$ satisfait la condition GFDA (pour {\em good finite dimensional approximation}) s'il existe un système compatible de projections $\pi_i :V\to W_i$ sur des espaces de dimension finie tel que
\begin{enumerate}
  \item $\n{v}=\lim_{i\to \infty}\n{\pi_i (v)}$.
  \item $\forall\epsilon>0$, pour tout suite $\rho_i$ décroissant vers 0, il existe $N$ tel que chaque fois que deux vecteurs $v$ et $v'$ satisfont, pour $i=1,\ldots,N$, $\n{v}-\n{\pi_i (v)}<\rho_i \n{v}$, $\n{v'}-\n{\pi_i (v')}<\rho_i \n{v'}$ et $\n{\pi_N (v)-\pi_N (v')}<N^{-1}\max\{\n{v},\n{v'}\}$, alors $\n{v-v'}<\epsilon\max\{\n{v},\n{v'}\}$.
\end{enumerate}
\end{defi}

\begin{exs}
\begin{enumerate}
  \item Cette classe contient $\ell^{1}$ mais pas $L^1$. En effet, GFDA $\Rightarrow$ RNP.
  \item Tout dual séparable possède une norme équivalente qui est GFDA. En effet, d'après Kadec et Klee, \cite{Kadec}, \cite{Klee}, dans un dual séparable, quitte à changer de norme, la convergence faible et la convergence des normes entraîne la convergence forte.
\end{enumerate}
\end{exs}

\begin{thm}
{\em (Cheeger-Kleiner, 2006, \cite{Cheeger-Kleiner-GFDA})}. Les espaces PI dont les cônes tangents ont une dimension de Hausdorff strictement supérieure à leur dimension topologique ne se plongent pas bi-lipschitz dans les espaces de Banach qui satisfont la condition GFDA. Noter qu'il existe des espaces PI qui se plongent dans $L^{1}$.
\end{thm}

\subsection{Le groupe d'Heisenberg}
\begin{defi}
Le {\em groupe d'Heisenberg} $\H$ est le groupe de Lie de dimension 3 dont l'algèbre de Lie a pour base $\xi$, $\eta$ et $\zeta$ tels que $[\xi,\eta]=\zeta$. Les champs de vecteurs invariants à gauche $\xi$ et $\eta$ engendrent un champ de plans $H$. La {\em distance de Carnot} $d(x,x')$ est l'$\inf$ des longueurs des courbes tangentes à $H$ joignant $x$ à $x'$. L'{\em homothétie} $\delta_{t}$ est l'automorphisme induit par $\delta_{t}(\xi)=t\xi$, $\delta_{t}(\eta)=t\eta$, $\delta_{t}(\zeta)=t^2 \zeta$. Elle multiplie les distances de Carnot par $t$.
\end{defi}

La finitude de la distance de Carnot résulte de la figure suivante.

\begin{center}
\includegraphics[width=1.5in]{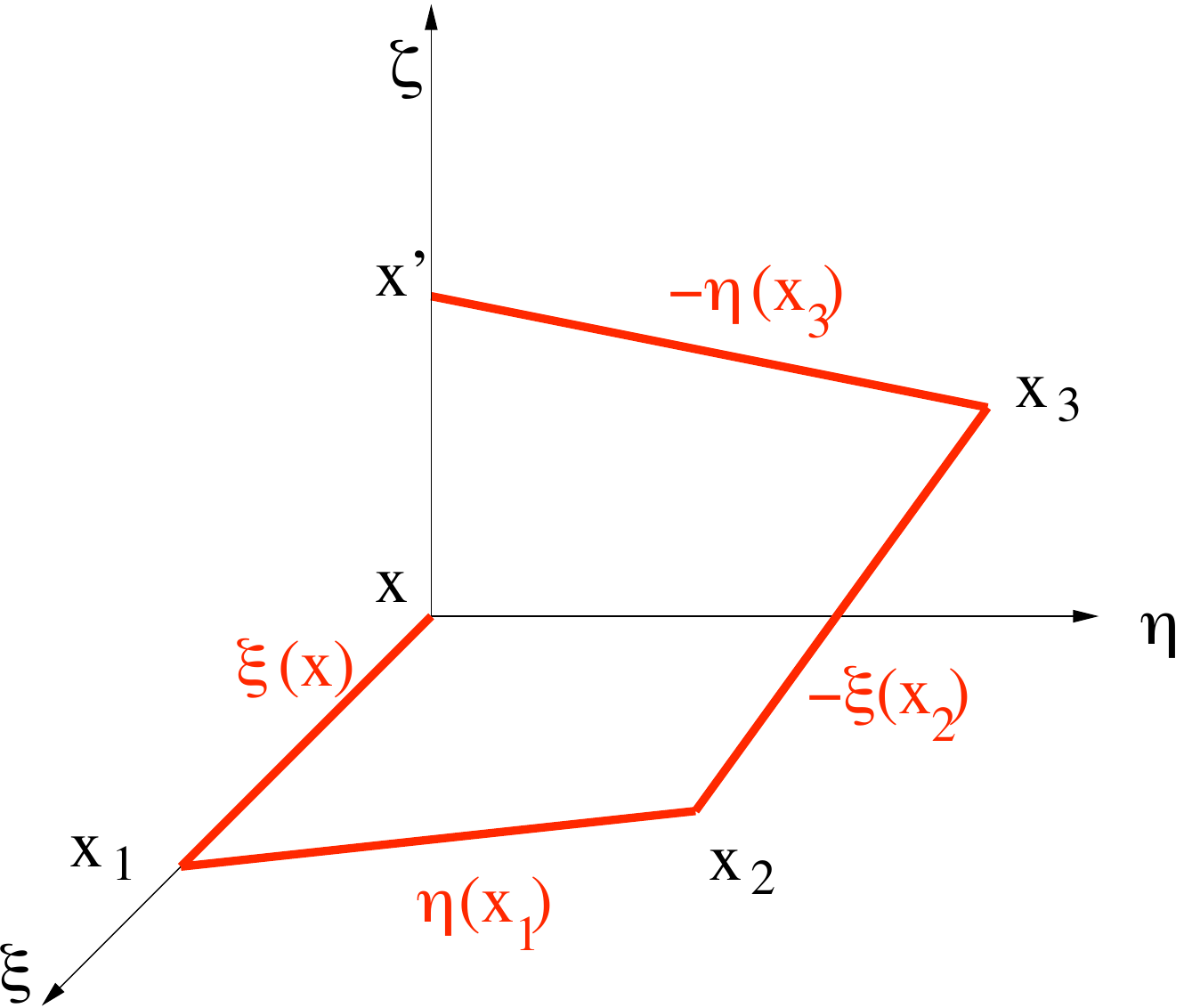}
\end{center}

\begin{rem}
\begin{enumerate}
  \item $d(x,x\exp(t^2 \zeta))=td(1,\exp(\zeta))=\mathrm{const.}\,t$.
  \item $\mathrm{volume} B(x,t)=t^{4}\mathrm{volume} B(x,1)=\mathrm{const.}\,t^4$, donc la dimension de Hausdorff est 4.
\end{enumerate}\end{rem}

\subsection{Non plongeabilité du groupe de Heisenberg dans les espaces de Banach RNP}

\begin{enumerate}
\item {\em Un application lipschitzienne possède des dérivées partielles horizontales presque partout}.

\smallskip

Par définition de la propriété de Radon-Nikodym.
\item {\em Ces dérivées sont approximativement continues presque partout}.

\smallskip

C'est général pour les espaces doublants.
\item {\em En un tel point $x$, $d(f(x),f(x'))=o(d(x,x'))$ si $x'$ est sur la droite verticale passant par $x$}.
\end{enumerate}

\smallskip

Soit $x'=x\exp(t^2 \zeta)$, de sorte que $d(x,x')\sim t$. On joint $x$ à $x'$ par des courbes intégrales de $\xi$, $\eta$, $-\xi$, $-\eta$, d'extrémités $x=x_0$, $x_1$, $x_2$, $x_3$, $x_4 =x'$ comme indiqué sur la figure ci-dessus. Alors 
$$f(x_1)-f(x)\sim t\xi f(x),\quad f(x_3)-f(x_2)\sim -t\xi f(x),$$ 
d'où
$$f(x_1)-f(x)+f(x_3)-f(x_2)=o(t).$$ 
De même $f(x_2)-f(x_1)+f(x')-f(x_3)=o(t)$. En additionnant, $f(x')-f(x)=o(t)$. Ceci est incompatible avec une inégalité de la forme $\parallel f(x')-f(x) \parallel \geq c\,d(x,x')$.

\section{Preuve de l'impossibilité de plonger dans $L^1$}

\begin{rem}
$L^1$ n'a pas la propriété de Radon-Nikodym. 
\end{rem}
En effet, $t\mapsto 1_{[0,t]}$, $\R_+ \to L^1(\R_+ )$, est isométrique, mais nulle part différentiable.

\subsection{Schéma de la preuve}

\begin{enumerate}
  \item A une application $f:X\to L^1 (Y,\nu)$, on associe une famille canonique de sous-ensembles $S\subset X$, généralisant les ensembles de niveau.
  \item $f$ est à variation bornée si et seulement si presque tout $S$ est de périmètre fini.
  \item {\em (Franchi, Serapioni, Serra-Cassano, 2001, \cite{Franchi-Serapioni-Serra-Cassano})}. Si $S\subset \H$ est de périmètre fini, alors, en presque tout point, les dilatés par $\delta_t$ de $S$ convergent vers des demi-espaces verticaux. 
  \item Au voisinage des bons points, à petite échelle, $f$ factorise approximativement par $\H/[\H,\H]$, ce qui l'empêche d'être bilipschitzienne. 
\end{enumerate}

\subsection{Distances de coupure}

\begin{defi}
Une {\em semi-distance de coupure élémentaire} sur $X$, c'est 
$$
d_S (x,x')=|1_{S}(x)-1_{S}(x')|
$$ 
où $S\subset X$ est une coupure. Une {\em semi-distance de coupure} est une somme de semi-distances de coupure élémentaires, i.e.
\begin{eqnarray*}
d(x,x')=\int_{\{S\}}d_S (x,x')\,d\mu_d (S)
\end{eqnarray*}
où $\mu_d $ est une mesure sur l'ensemble des coupures.
\end{defi}

Une référence générale sur les distances de coupure est \cite{Deza-Laurent}.

\begin{lem}
\label{a}
{\em (P. Assouad, 1977, \cite{Assouad-77})}. Une semi-distance $d$ sur $X$ est induite par une application $f:X\to L^1 (Y,\nu)$ si et seulement si c'est une semi-distance de coupure.
\end{lem}

\proof ~\textbf{du lemme d'Assouad}, $\Leftarrow$. Supposons que $d$ est une semi-distance de coupure. Fixons une origine $o\in X$. Soit $S(x)$ l'ensemble des coupures qui séparent $x$ de $o$. Alors $x\mapsto 1_{S(x)}$ plonge $(X,d)$ isométriquement dans $L^1 (\{S\},\mu_d )$.

\subsection{La mesure de coupure}

\begin{defi}
L'{\em épigraphe} d'une fonction $u:Y\to\R$ est 
$$
E_u =\{(y,t)\in Y\times\R\,|\,t^{-1}u(y)>1\}.
$$
\end{defi}

\begin{center}
\includegraphics[width=1.8in]{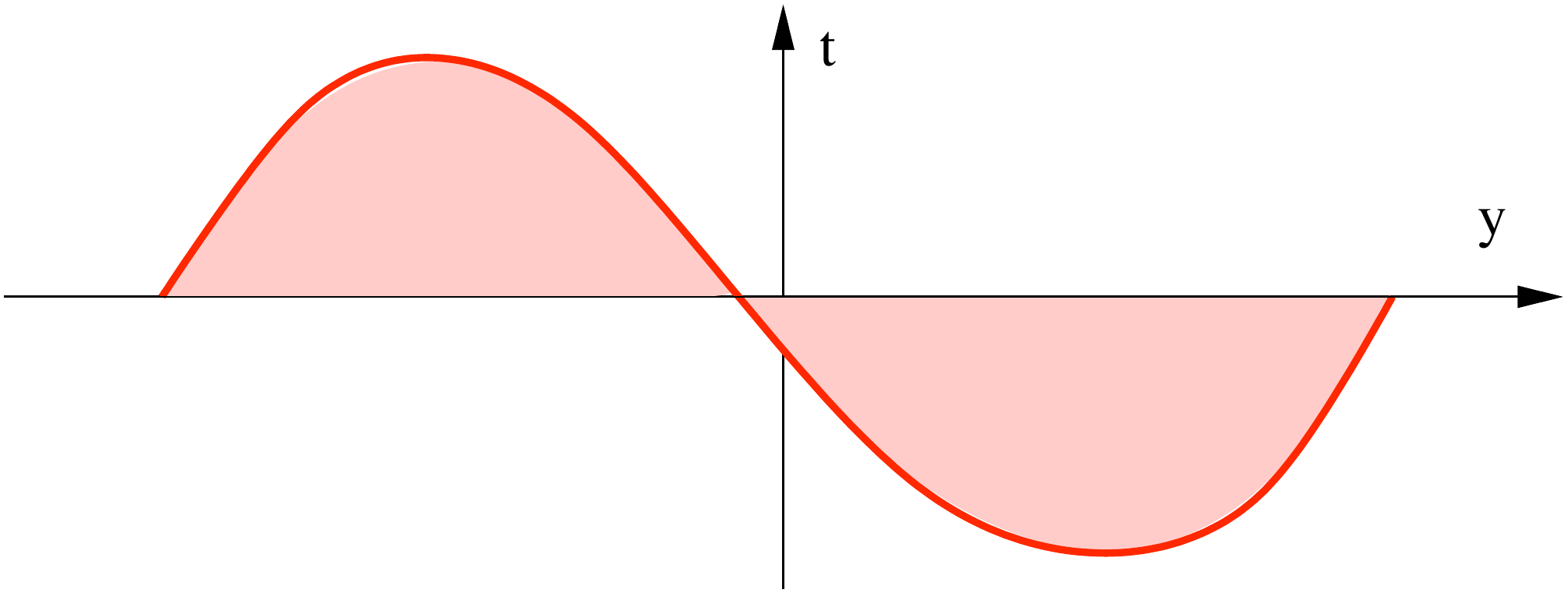}\hskip1cm\includegraphics[width=1.8in]{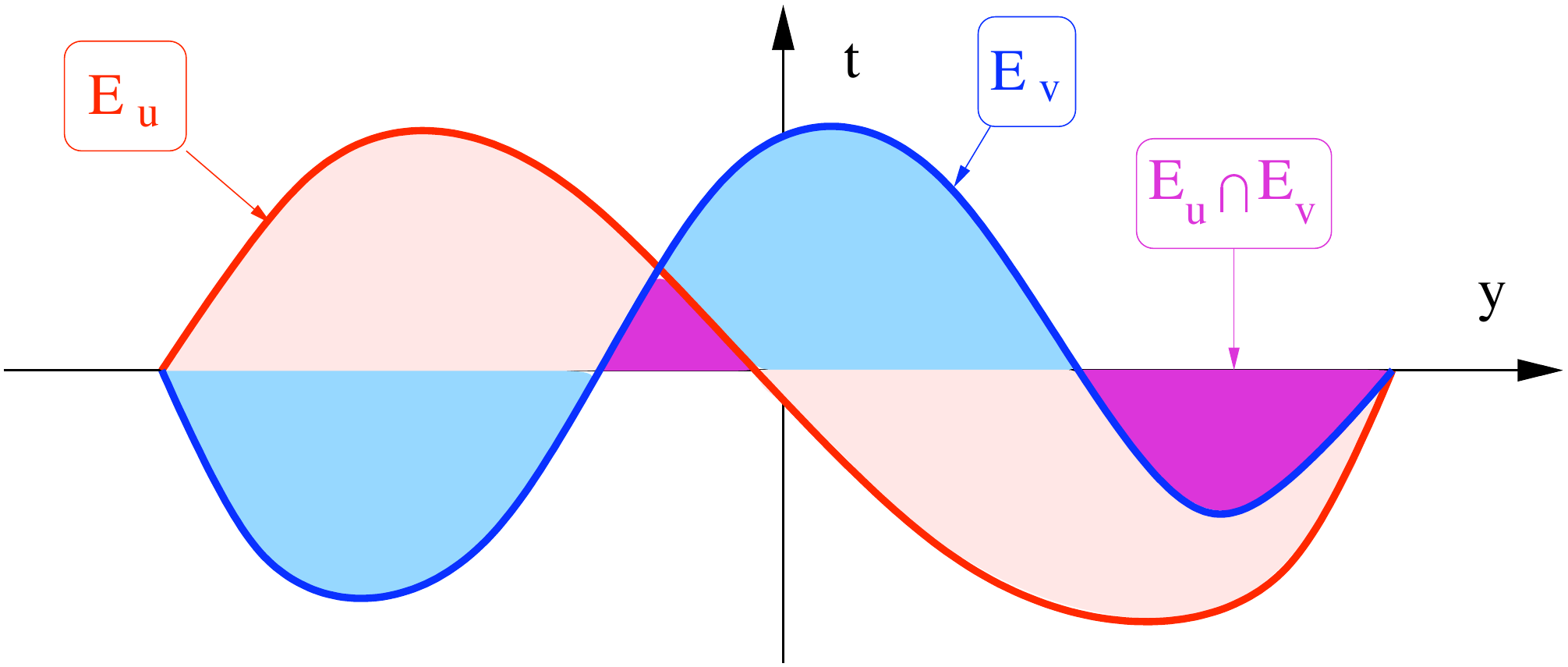}
\end{center}

\begin{lem}
Si $u$, $v\in L^1 (Y,\nu)$, $\n{u-v}_{L^1}=(\nu\otimes dt)(E_{u}\Delta E_{v})$.
\end{lem}

\proof ~\textbf{du lemme d'Assouad}, $\Rightarrow$. Soit $f:X\to L^1 (Y,\nu)$. A chaque $(y,t)\in Y\times\R$, correspond une coupure $S(y,t)=\{x\in X\,|\,(y,t)\in E_{f(x)}\}=\{x\in X\,|\,t^{-1}f(x)(y)>1\}$. Soit $\mu_{f}=S_* (\nu\otimes dt)$. Alors pour tous $x$, $x'\in X$,
\begin{eqnarray*}
d(f(x),f(x'))&=&\n{f(x)-f(x')}_{L^1}=(\nu\otimes dt)(E_{f(x)}\Delta E_{f(x')})\\
&=&\int_{Y\times\R}|1_{E_{f(x)}}(y,t)-1_{E_{f(x')}}(y,t)|\,d\nu(y)\,dt\\
&=&\int_{Y\times\R}|1_{S(y,t)}(x)-1_{S(y,t)}(x')|\,d\nu(y)\,dt\\
&=&\int_{\{S\}}|1_{S}(x)-1_{S}(x')|\,d\mu_{f}(S).
\end{eqnarray*}

\subsection{Variation bornée}

On rappelle que si $h$ est une fonction lisse sur un ouvert $U$ de $\R^n$,
$$
\mathrm{aire}(graphe(h))=\int\sqrt{1+|\nabla h|^2}\sim Vol(U)+\n{\mathrm{Lip}h}_{L^1}. 
$$
Toute limite $L^1$ de fonctions $h$ de norme $\n{\mathrm{Lip}h}_{L^1}$ bornée conduit à une surface d'aire finie. Cela a conduit Tonelli à la
\begin{defi}
Soient $X$ et $Y$ des espaces métriques. Une application $L^1$ $f:X\to Y$ est {\em à variation bornée} si c'est la limite $L^1$ d'une suite de fonctions localement lipschitziennnes $h_i$ telles que $\mathrm{Lip}h_i$ reste borné dans $L^1$. La borne inférieures des limites des $\int \mathrm{Lip}h_i$ sur toutes les approximations $h_i$ s'appelle la {\em variation} de $f$.

\medskip

Un sous-ensemble $S\subset X$ a un {\em périmètre fini} si sa fonction caractéristique $1_S$ est à variation bornée. Le {\em périmètre} de $S$ est égal à la variation of $1_S$.
\end{defi}

Pour les fonctions $BV$ à valeurs réelles, la {\em formule de la coaire} donne
\begin{eqnarray*}
\mathrm{variation}(h)=\int_{\R}\textrm{périmètre}(\{h>t\})\,dt.
\end{eqnarray*}
Cela s'étend aux applications $X\to L^1 (Y,\nu)$ comme suit.

\begin{thm}
{\em (Cheeger-Kleiner, \cite{Cheeger-Kleiner-L1})}. Soit $X$  un espace PI. Soit $f\in L^{1}(X,L^{1}(Y,\nu))$. Alors $f$ est à variation bornée si et seulement si $\mu_f$-presque toute coupure est de périmètre fini. De plus
\begin{eqnarray*}
\int_{\{S\}}\textrm{périmètre}(S)\,d\mu_f (S)=\int_{Y}\mathrm{variation}(f(\cdot,y))\,d\nu(y)\leq \mathrm{const.~variation}(f).
\end{eqnarray*}
\end{thm}

\subsection{Mauvais points}

Désormais, $X=\R^n$ ou $X=\H$. Un {\em demi-espace} dans $\H$ est l'image réciproque de $\R_+$ par un homomorphisme de groupes $\H\to\R$ (il est borné par un plan \textbf{vertical}). 

\begin{nota}
Pour $S\subset X$, $x\in\partial S$, soit 
$$
\alpha(S,x,r)=\min_{H\mathrm{\, demi-espace\, passant\, par\, }x}~\oint_{B(x,r)}|1_S -1_H|.
$$
On dit qu'un point $x$ est $(\epsilon,R)$-mauvais pour $S$, et on note $x\in Bad_{\epsilon,R}(S)$, s'il existe $r<R$ tel que $\alpha(S,x,r)>\epsilon$.
\end{nota}

Pour tout $S$, le {\em périmètre} de $S$ est une mesure $\n{\partial S}$ sur $X$, portée par la frontière de $S$. On peut la restreindre aux mauvais points. Etant donnée une mesure $\mu$ sur l'ensemble des coupures de périmètre fini, la {\em mesure du mauvais périmètre} $\lambda_{\mu,\epsilon,R}$ est définie sur une fonction continue $u$ par
\begin{eqnarray*}
\int_{X}u(x)\,d\lambda_{\mu,\epsilon,R}(x)=\int_{\{S\}}\int_{Bad_{\epsilon,R}(S)}u(x)\,d\n{\partial S}(x)\,d\mu(S).
\end{eqnarray*}

\begin{thm}
\label{i}
{\em (Franchi, Serapioni, Serra-Cassano, 2001, \cite{Franchi-Serapioni-Serra-Cassano})}.
Soit $S$ un ensemble de périmètre fini dans $\H$. Le périmètre de l'ensemble $Bad_{\epsilon,R}(S)$ tend vers 0 quand $R$ tend to 0.
\end{thm}

\begin{cor}
Soit $\mu$ une mesure sur l'ensemble des coupures de périmètre fini. La masse totale de la mesure $\lambda_{\mu,\epsilon,R}$ tend vers 0 quand $R$ tend to 0.
\end{cor}

\subsection{Approximation d'une mesure de coupure par une mesure portée par les demi-plans}

Soit $\mu$ une mesure sur l'ensemble des coupures de périmètre fini, soit $d_\mu$ la semi-distance correspondante. Etant donnés $x\in X$ et $r>0$, soit $\delta_{x,r}^* d_\mu$ la distance distance composée avec l'homothétie de rapport $r$ et de centre $x$.

\begin{thm}
{\em (Cheeger-Kleiner)}. Pour presque tout $x\in X$, il existe des mesures $\mu_r$ portées par les demi-espaces telles que $\n{\frac{1}{r}\delta_{x,r}^* d_\mu -d_{\mu_{r}}}$ tend vers 0 dans $L^1 (X\times X)$.
\end{thm}

\proof. Par différentiation de la mesure du mauvais périmètre, on obtient pour presque tout $x$ un ensemble de mesure presque pleine dans la boule $B(x,r)$ de points où la plupart (au sens de $\mu$) des coupures sont proches de demi-espaces. Pour une telle coupure $S$, on choisit le demi-espace $HS(S)$ le plus proche de $S$, et on pose $\mu_r =\frac{1}{r}(HS\circ\delta_{x,r})_* \mu$.

\proof~\textbf{du théorème de non-plongement}. Si $f:\H\to L^1$ est un plongement bilipschitzien, de mesure de coupure $\mu$, $d_\mu (x',x'')=d(f(x'),f(x''))\geq \mathrm{const.}d(x',x'')$, donc 
$$\frac{1}{r}\delta_{x,r}^* d_\mu (x',x'')\geq \mathrm{const.}d(x',x'').$$ 
En revanche, une distance de coupure concentrée sur les demi-espaces satisfait
$$d_{\mu_{r}}(x',x'')=d_{\mu_{r}}(x'~\mathrm{mod}~Z(\H),x''~\mathrm{mod}~Z(\H)).$$ 
Deux telles semi-distances ne peuvent pas être $L^1$-proches.

\subsection{La normale unitaire d'un ensemble de périmètre fini}

On donne la preuve du théorème \ref{i}. Ici, $X=\R^n$ or $\H$. Soit $S\subset X$ un ensemble de périmètre fini. On utilise la formule de la divergence pour définir la normale au bord. Dans $\H$, la divergence d'un champ de vecteurs horizontal est définie comme suit : $\phi=\phi_\xi \xi +\phi_\eta \eta$,
$div(\phi)=\xi\phi_\xi +\eta\phi_\eta$.

\begin{lem}
{\em (De Giorgi, 1954, \cite{DeGiorgi})}. Il existe un champ de vecteurs unitaire (horizontal) mesurable $\nu$ tel que pour tout champ de vecteur à support compact (horizontal) $\phi$,
\begin{eqnarray*}
-\int_{S}div(\phi)=\int_{X}\langle\nu,\phi\rangle d\n{\partial S}.
\end{eqnarray*}
\end{lem}

\proof. Pour une fonction lipschitzienne $h$, la formule de la divergence entraîne
\begin{eqnarray*}
|\int_{S}h\,div(\phi)|\leq\n{\phi}_{L^{\infty}}\int_{X}\mathrm{Lip}h.
\end{eqnarray*}
Cela montre que $\phi\mapsto -\int_{S}div(\phi)$ est une mesure de Radon à valeurs vecteorielles dont la masse est bornée par la variation de $(h)$. Ce fait s'étend aux fonction à variation bornée, en particulier à $1_S$. Dans ce cas, la mesure de Radon est absolument continue par rapport à la mesure de périmètre, et $\nu$ est sa densité (Riesz).

\subsection{Rectifiabilité des ensembles de périmètre fini}

\begin{lem}
{\em (Ambrosio, 2001, \cite{Ambrosio})}. En $\n{\partial S}$-presque tout point $x$, la $\n{\partial S}$-mesure d'une boule $B(x,r)$ est $\sim r^3$.
\end{lem}
Cela entraîne que les fonctions mesurables comme $\nu$ sont approximativement continues presque partout.

\medskip

Par compacité, on peut extraire une sous-suite convergente des dilatés $\delta_{x,1/r}(S)$. La limite est une ensemble $E$ de périmètre localement fini. Les normales unitaires convergent aussi $\n{\partial S}$-presque partout, donc la normale unitaire de $E$ est presque partout constante.

\begin{lem}
Un ensemble $E$ de périmètre localement fini dont la normale unitaire est presque partout égale à $\xi$ est un demi-espace (vertical).
\end{lem}
En effet, si on se déplace depuis l'origine, dans les deux sens le long des orbites de $\eta$, et positivement le long des orbites de $\xi$, on atteint exactement tous les points d'un demi-espace (vertical).  

Cela complète en gros la preuve du théorème \ref{i}. Franchi, Serapioni and Serra-Cassano vont plus loin dans la théorie de la rectifiabilité.

\begin{thm}
{\em (Franchi, Serapioni et Serra-Cassano, 2001, \cite{Franchi-Serapioni-Serra-Cassano})}. A un ensemble de mesure de Hausdorff 3-dimensionnelle nulle près, le bord d'un ensemble de périmètre fini est une réunion dénombrable de morceaux compacts de surfaces définies par des équations $g=0$ de classe $C^1$, dont le gradient horizontal ne s'annule pas.
\end{thm}

\section*{Appendice : Réduction du non plongement des boules au non plongement de $(\H,d)$}

\begin{prop}
L'impossibilité de plonger le groupe d'Heisenberg dans $L^1$ entraîne que $L_n$ tend vers $+\infty$.
\end{prop}

\proof. 
Le lemme suivant remonte à S. Kakutani (et résulte du lemme \ref{a}).

\begin{lem}
{\em (S. Kakutani, 1939, \cite{Kakutani})}. L'ensemble des semi-distances plongeables dans $L^1$ est fermé. De plus, un espace semi-métrique $X$ se plonge isométriquement dans $L^1$ si et seulement si tout sous-ensemble fini de $X$ se plonge isométriquement dans $L^1$.
\end{lem}

Supposons que la boule de rayon $n$ de $(\H_{\Z},d_{\Z})$ admette un plongement $L$-bilipschitzien dans $L^1$. Comme les distances induites $d_n$ sont bornées, on peut supposer qu'elles convergent simplement vers une distance $d'$. Alors $d'$ se plonge dans $L^1$, elle est $L$-équivalente à $d_{\Z}$. Donc $(\H_{\Z},d_{\Z})$ admet un plongement $L$-bilipschitzien dans $L^1$. 

\begin{lem}
{\em (S. Kakutani, 1939, Bretagnolle, Dacunha-Castelle, Krivine, 1966, \cite{Bretagnolle-Dacunha-Castelle-Krivine})}.
Tout ultraproduit d'espaces $L^p$ est à nouveau un espace $L^p$.
\end{lem}

En particulier, tout cône asymptotique de $L^1$ est isométrique à $L^1$. Par conséquent, un plongement bilipschitzien (ou quasiisométrique) de $(\H_{\Z},d_{\Z})$ dans $L^1$ donne un plongement bilipschitzien de $(\H,d)$ dans $L^1$, contradiction.

\vskip1cm

\noindent  Pierre Pansu\\ Laboratoire de Math{\'e}matiques d'Orsay\\
UMR 8628 du CNRS\\
Universit{\'e} Paris-Sud\\
91405 Orsay C\'edex, France\\
\texttt{pierre.pansu\char`\@ math.u-psud.fr}\\
\texttt{http://www.math.u-psud.fr/$\sim$pansu}
  
\end{document}